\documentclass[11pt,reqno]{article}
\usepackage{amsmath,amsthm,amsfonts,amssymb,amscd}

\def\got k{\operatorname{{\cal K}}}

\newtheorem{Theorem}{Theorem}
\newtheorem{Corollary}{Corollary}
\newtheorem{Proposition}{Proposition}
\newtheorem{Lemma}{Lemma}
\newtheorem{Claim}{Claim}
\theoremstyle{Definition}
\newtheorem{Definition}{Definition}
\theoremstyle{Definition}
\newtheorem{Example}{Example}
\newtheorem{Question}{Question}

\theoremstyle{Remark}

\theoremstyle{Remark}
\newtheorem{Remark}{Remark}
\def\fa{{\mathcal{F}}}
\def\Aff{\operatorname{{Aff}}}
\def\Cod{\operatorname{{cod}}}
\def\sing{\operatorname{{Sing}}}

\def\dif{\operatorname{{Diff}}}

\def\Hol{\operatorname{{Hol}}}
\def\O{{\mathcal{O}}}
\def\C{{{\mathbb C}}}
\def\codim{\operatorname{{codim}}}

\title{Foliations invariant under Lie group transverse actions}

\author{Alexandre Behague and Bruno Sc\'ardua}

\date{}
\begin{document}
\maketitle

\begin{abstract}
In this paper we study  (smooth and holomorphic) foliations which
are invariant under transverse actions of Lie groups.
\end{abstract}


\section{Introduction and main results}
\label{section:Introduction}

In the study of foliations it is very useful to consider the
transverse structure\footnote{We are in debt with Professor J. J.
Duistermaat for  various suggestions and valuable
remarks.}\footnote{MSC Classification: 57R30, 22E15,
22E60.}\footnote{Keywords: Foliation, Lie transverse structure,
fibration.}. Among the simplest transverse structures are Lie
group transverse structure, homogeneous transverse structure and
Riemannian transverse structure. In the present work we consider a
sightly different situation; foliations which are invariant under
Lie group transverse actions. Another motivation for this work is
the well-known result of Tischler \cite{Tischler} asserting that
if  a closed oriented manifold admits
 a (codimension one) foliation  which is invariant under a transverse flow
 then the manifold is a fiber bundle over the circle.
 In this work we look for generalizations of these result for higher codimension
foliations. All manifolds are assumed to be connected and
 oriented. All foliations are assume to be smooth, oriented and transversely oriented.

 Let $M$ be a manifold, $\fa$ a
codimension $q$ foliation on $M$ and $G$ a Lie group of dimension
$\dim G=\codim \fa=q$. We shall also say that $\fa$ is  {\it
invariant under a transverse action} of the group $G$, $\fa$ is
$G$-i.u.t.a. for short, if there is an  action $\Phi\colon G\times
M\rightarrow M$ of $G$  on $M$ such that: (i) the action is {\it
transverse to $\fa$}, i.e.,   the orbits of this action have
dimension $q$ and intersect transversely the leaves of $\fa$ and
(ii) $\Phi$ {\it leaves   $\fa$ is  invariant}, i.e.,  the maps
$\Phi_{g}: x \mapsto \Phi(g,x)$  take leaves of $\fa$ onto leaves
of $\fa$.

Let $\fa$ be a foliation on $M$ such that $\fa$ is $G$-i.u.t.a. It
is not difficult to prove the existence of a {\it Lie foliation
 structure  for $\fa$ on $M$ of model $G$} in the sense of
Ch. III, Sec. 2 of \cite{Godbillon}. We shall then say that $\fa$
{\it has $G$-transversal structure} and prove (with a
self-contained proof) the existence of a {\it development} for
$\fa$ as in Proposition 2.3, page 153 of \cite{Godbillon} (Ch.III,
Sec. 2).  Indeed, we have a sort of strong form of this procedure
in Section~\ref{Section:development} with a self-contained proof
(Proposition~\ref{PropositionDevelopment})

Indeed, from the proof of Proposition~\ref{PropositionDevelopment}
we obtain an {\it algebraic model} for any foliated manifold
$(M,\fa)$ assuming that $\fa$ is $G$-i.u.t.a. Given a leaf $L$ of
$\fa$ we define $H(L)$ as the set of $g\in G$ such that
$\Phi(g,l)\in L$ for every {\rm(}or equivalently for some{\rm)}
$l\in L$. Then $H(L)$ is a (not necessarily closed) subgroup of
$G$ which we provide  the discrete topology. We have the following
algebraic model for the general foliation invariant under a
transverse Lie group action.
\begin{Theorem}[Algebraic model]
\label{Theorem:algebraicmodel} Let $\fa$ be a foliation and
$G$-i.u.t.a. Given a leaf $L$ of $\fa$ there is a natural proper
free  action  of $H(L)$ on $G\times L$ with a smooth  quotient
manifold $(G\times L)/H(L)$, which is $G$-equivariantly
diffeomorphic to $M$.  The leaves of $\fa$  are the sets
$\Phi(p(\{g\}\times L))$, $g\in G$ where $p:G\times L\to (G\times
L)/H(L)$ denotes the canonical  projection.
\end{Theorem}

As a consequence of the above construction we have:

\begin{Theorem}[Fibration theorem] \label{Theoremfibrebundlenew}
Let $M$ be a connected manifold and ${\cal F}$ a foliation in $M$
which is invariant under a transverse action of a Lie group $G$.
Then the following statements are equivalent:

\noindent{\rm (a)} ${\cal F}$ has a leaf $L$ which is closed in
$M$.

\noindent{\rm (b)} $H(L)$ is a discrete {\rm(}i.e., closed and
zero-dimensional{\rm)} subgroup of $G$.

\noindent {\rm (c)} The projection $\pi:G\times L\to G$ onto the
first factor induces a smooth fibration $M\simeq G\times_H(L) L\to
G/H(L)$, of which the fibers are the leaves of ${\cal
F}$.\end{Theorem}

From Theorem~\ref{Theoremfibrebundlenew} we immediately obtain:

\begin{Corollary} If $L$ is compact, then $H(L)$ is finite, and we have a fiber
bundle over $G/H(L)$.
\end{Corollary}

Additional consequences of Theorem~\ref{Theoremfibrebundlenew}
are:

\begin{Corollary}
\label{Corollary:finitehomotopy} If $\pi_1(M,x)$ is finite, then
$H(L)$ is closed and  $\fa$ is a fibration over the base space
$G/H(L)$. If moreover $M$ is compact, then  $L$ and  $G$ are
compact. In case $G$ is simply connected the latter also implies
that $G$ is semi-simple.
\end{Corollary}

\begin{Corollary}\label{Corollary:codimensiontwo}
Let $M$ be a compact manifold supporting a
codimension two foliation $\fa$ invariant under a Lie group
transverse action. If $\fa$ has a compact leaf then $M$ is a fibre
bundle over the two torus.
\end{Corollary}

A well-known consequence of (the proof of) Tischler's theorem is
that {\em on a  compact, connected, oriented manifold $M$ with
$\dim H_1(M, R)\leq 1$ or, equivalently, with $\dim H^1_{de
Rham}(M)\leq 1$, any foliation  ${\cal F}$ of a codimension one
foliation which admits an $\mathbb R$-transversal structure is
defined by a fibration over a circle and in particular $\dim
H^1(M)=1$}. This   can be generalized as follows:

\begin{Corollary}
\label{Corollary:d-torus} Let $M$ be a compact smooth manifold and
${\cal F}$ a foliation in $M$ which admits a $G$-transversal
structure. If $d:=\dim H^1_{de Rham}(M)\leq \dim {{\cal L}G} -
\dim D{{\cal L}G}$, then we have equality  and $M$ fibers over a
$d$-dimensional torus in such a way that the leaves of ${\cal F}$
are contained in the fibers of this fibration. In particular, if
${{\cal L}G}$ is abelian and $\dim H^1_{de Rham}(M)\leq \dim
{{\cal L}G}$, then ${\cal F}$ is a fibration over a torus.
\end{Corollary}

Regarding (codimension two) foliations which are  $\Aff({\mathbb
R})$-i.u.t.a. it is not difficult to prove that if \,
$H^{1}(M,{\mathbb R})=0$ then $\fa$ is given by a submersion
$F\colon M\to \Aff({\mathbb R})$. This fact admits the following
strong form:

\begin{Theorem}
\label{Theorem:solvable} If $M$ is connected, ${\cal F}$ is a
foliation which is invariant under the transverse action of a
simply connected {\underline{solvable}} Lie group $G$, and $H^1(M,
\mathbb R)=0$, then $M$ is diffeomorphic to $G\times L$, where $L$
is any leaf of $\fa$, and the foliation is given by the projection
to the first factor $G$, which is diffeomorphic to a vector space.
\end{Theorem}

\subsubsection*{Holomorphic foliations}

 In Section~\ref{section:Holomorphic} we carry out  the study of
 the holomorphic case and prove
 an analogous to  Theorem~\ref{Theoremfibrebundlenew}
 (Theorem~\ref{Theorem:locallyholomorphictrivialfibration}).
 For  the case $M$ is a compact K\"ahler manifold we prove
 that  if $\fa$ has a $G$-transverse structure then
 the universal covering of $G$ is isomorphic to $(\C^{q},+)$
(Proposition~\ref{Proposition:KahlerLiestructure}), if moreover
$\fa$ has some compact leaf then  $G=\C^{q}/H$ for some closed
subgroup $H < \C^{q}$
(Proposition~\ref{Proposition:complexcovering}).  Finally, we
consider a codimension one algebraic foliation ${\fa}_{0}$ on
$\C^{n}$, we denote by $\fa$ its extension to $\C\mathbb P^{n}$.
We prove an  following extension result
(Theorem~\ref{Theorem:extension}) which implies the following:

\begin{Theorem}\label{Theorem:extensionaction}Let $\fa$ be a
codimension $q$ singular holomorphic foliation on $\C\mathbb
P^{n}$ and suppose that there is  an algebraic irreducible
hypersurface $\Gamma \subset \C\mathbb P^{n}$ which is not
$\fa$-invariant  and a holomorphic action of a Lie group $G$ on
$\C\mathbb P^{n}\setminus \Gamma$ transverse to $\fa$ and under
which $\fa$ is invariant.  Then the action  extends to an action
on $\C\mathbb P^{n}$ and, in particular, $G$ embeds as a
{\rm(}linear{\rm)} subgroup of birational maps of $\C\mathbb
P^{n}$.
\end{Theorem}

\noindent{\bf Acknowledgement}: This paper is based on an original
manuscript of the first named author on ``Foliations invariant
under Lie group transverse actions" and on a number of mails from
Professor J.J. Duistermaat. We are very grateful to Professor J.J.
Duistermaat for reading the original manuscript and for suggesting
and sketching various improvements on the original results. In
particular, the construction of the {\it algebraic model} in
Section~\ref{Section:algebraicmodel} and the strong forms
Theorems~\ref{Theorem:solvable} and
Corollary~\ref{Corollary:d-torus} are due to him. Also due to him
is the introduction and study of the notion of invariance under a
local action in Section~\ref{Section:localaction}.

\section{Examples}


In Section~\ref{Section:algebraicmodel} we construct the algebraic
model of the general foliation invariant under a Lie group
transverse action. This provides a number of examples of foliated
spaces with invariant foliations. Below we give some more concrete
examples:

\begin{Example}{\rm The most trivial example of a foliation
invariant a Lie group transverse action is given by the product
foliation on a manifold $M=G\times N$ product of a Lie group $G$
by a manifold $N$. The leaves of the foliation are of the form
$\{g\}\times N$ where $g\in G$.}
\end{Example}

\begin{Example}{\rm Let $H$ be a closed
(normal) subgroup of a Lie group $G$. We consider the action
$\Phi\colon H\times G\rightarrow G$ given by $\Phi(h,g)=h.g$ and
the quotient map $\pi\colon G\rightarrow G/H$ (a fibration) which
defines a foliation $\fa$ on $G$. Given
$x\in{\fa}_{g}=\pi^{-1}(Hg)$ we have $\pi(x)=Hg$ and
$\Phi_{h}(x)=h.x$. But $\pi(\Phi_{h}(x))=\pi(h.x)=H.hx=Hx$ implies
that $\Phi_{h}(x)\in \pi^{-1}(Hx)={\fa}_{x}$ and the orbit
$\O(g)=Hg$ is transverse to the fiber $\pi^{-1}(Hg)$. Hence, $\fa$
is a foliation invariant under the transverse action $\Phi$. Now
let $G$ be a simply-connected group, $H$ a discrete subgroup of
$G$ and $\phi\colon H\rightarrow \dif(G)$ the natural
representation given by $\phi(h)=L_{h}$. The universal covering of
$G/H$ is $G$ with projection $\pi\colon G\rightarrow G/H$ and we
have $\pi_{1}(G/H)\simeq H$ because $\pi \circ f(g)=Hf(g)=Hg$ for
$f\in Aut(G)$, so $f(g)\simeq g$ implies that $f(g).g^{-1}\in H$
and $f(g)=h.g$, for some unique $h\in H$. Therefore $f=L_{h}$ and
then we define $Aut(G)\rightarrow H;\ f\mapsto h$, which is an
isomorphism. So, we may write $\phi\colon \pi_{1}(G/H)\rightarrow
\dif(G)$ and $\Phi\colon \pi_{1}(G/H)\times G\rightarrow G$. The
map $\Psi\colon H\times G\times G\rightarrow G\times G$ given by
$\Psi(h,g_{1},g_{2})=(L_{h}(g_{1}),L_{h}(g_{2}))$ is a properly
discontinuous action and defines a quotient manifold $M={G\times
G\over \Psi}$, which equivalence classes are the orbits of $\Psi$.
We have the following facts: (1) There exists a fibration
$\sigma\colon M\rightarrow G/H$ with fiber $G$, induced by
$\pi\colon G\rightarrow G/H$, and structural group isomorphic to
$\phi(H)<\dif(G)$. (2) The natural foliation $\fa$ on $G$ given by
classes $Hg;\ g\in G,$ is $\Phi$-invariant, such that the product
foliation $G\times \fa$ on $G\times G$ is $\Psi$-invariant and
induces a foliation ${\fa}_{0}$ on $M$, called suspension of $\fa$
for $\Phi$, transverse to $\sigma\colon M\rightarrow G/H$.}
\end{Example}

\begin{Example}{\rm Let $G={\mathbb P}SL(2,{\mathbb R})$ and $H={\Aff}({\mathbb R})\triangleleft G$. An
element of $G$ has the expression $x\mapsto {ax+b\over
cx+d}={a+{b\over x}\over c+{d\over x}}$. The group $H$ is the
isotropy group of $\infty$, so ${a\over c}=\infty \Leftrightarrow
c=0$ and an element of $H$ is given by $x\mapsto {ax+b\over
d}\simeq \left(\begin{array}{cc}
  a & b \\
  0 & {1\over a} \\
\end{array}\right)$. Since $H\triangleleft G,\ G$ has dimension 3 and $H$ has
dimension 2, we conclude that $G/H$ has dimension one. Thus we
have a fibration ${\mathbb P}SL(2,{\mathbb R})\rightarrow {\mathbb
R}P(1)\simeq S^{1}$ which is invariant under an action of
${\Aff}^{+}({\mathbb R})$ on ${\mathbb P}SL(2,{\mathbb R})$ having
leaves diffeomorphic to ${\mathbb R}^{+}\times {\mathbb R}\simeq
{\mathbb R}^{2}$.}
\end{Example}

\section{Foliations with Lie transverse structure}
Throughout this paper, except if explicitly mentioned otherwise,
$\mathbb F$ will denote the tangent bundle of the foliation $\fa$.
It is therefore an integrable subbundle of the tangent bundle $TM$
of $M$ and its connection form is flat, because of the
integrability.

\begin{Definition}[\cite{Godbillon}, Ch.III, p. 152]
\label{Definition:Liestructure} {\rm Let $G$ be a dimension $q$
Lie group and $\fa$ a codimension $q$ foliation on a
differentiable manifold $M$. A {\it Lie transverse structure of
model $G$ for $\fa$} is given by: {\rm(1)} An open cover $\{
U_{i}\}_{i\in I}$ of $M$ and a family of submersions $f_{i}\colon
U_{i}\rightarrow G$ such that ${\fa}|_{U_i}$ is given by
$f_{i}=$constant and, {\rm(2)} a family of locally constant maps
$\gamma_{ij}\colon U_{i}\cap U_{j}\rightarrow \{ \mbox{left
translations on }G\}$ such that
$f_{i}(x)=\gamma_{ij}(x).f_{j}(x),\ \forall x\in U_{i}\cap
U_{j}$.}
\end{Definition}

According to Ch.III, Cor. 2.4 and Prop. 2.7 in \cite{Godbillon},
the existence of a $G$-transversal structure for  $\fa$  is
equivalent to the existence of a ${\cal L}G$-valued smooth
differential one-form $\omega$ on $M$, such that the tangent
bundle $\mathbb {F}$ of ${\cal F}$ is equal to the kernel of
$\omega$, and $(d\omega)(u,v)= -[\omega(u),\omega(v)]$ for every
pair of vector fields $u$, $v$. Here ${\cal L}G$ denotes the Lie
algebra of the Lie group $G$, and, for every $x\in M$, $\mathbb
{F}_x=T_x{\cal F}_x$, the tangent space at $x$ of the leaf passing
through the point $x$. After the choice of a basis in ${{\cal
L}G}$, this amounts to having the suitable systems of differential
one-forms as follows:

\begin{Proposition}\label{Proposition:Lieforms} Let $M$ be a manifold
equipped with a  codimension $q$ having transverse structure of
model $G$. Then there exists an integrable system
$\{\Omega_{1},...,\Omega_{q}\}$ of one-forms defining $\fa$ on $M$
with $d\Omega_{k}=\Sigma_{i<j}c_{ij}^{k}\Omega_{i}\wedge
\Omega_{j}$, where $\{c_{ij}^{k}\}$ are the  structure constants
of the Lie algebra of $G$ for a certain basis.
\end{Proposition}

Lie foliations  exhibit the following structure (cf.
\cite{Godbillon} Prop. 2.9, p.155):

\begin{Proposition}[\cite{Godbillon}]\label{Proposition:Development}
For a foliation $\fa$ on $M$ the existing Lie transverse
structures of model $G$ are classified by: {\rm(i)} A Galoisian
covering map $\pi\colon P\rightarrow M$, {\rm(ii)} a homomorphism
$h\colon \pi_{1}(M)\rightarrow G$ such that $\pi^{-1}(e_G)\simeq
Aut(\pi)$ and {\rm(iii)} a submersion $\Theta\colon P\rightarrow
G$ which is a first integral for the pull-back foliation
$\pi^{*}{\fa}$ and is equivariant by $h$, that is,
$\Theta(\alpha.x)=h(\alpha)\Theta(x),\ \forall x\in M,\ \alpha \in
\pi_{1}(M)$. We call $(P,h,\Theta)$ a development of $\fa$. Two
developments $(P_{1},h_{1},\Theta_{1})$ and
$(P_{2},h_{2},\Theta_{2})$ define the same Lie group transverse
structure  if and only if there  is a diffeomorphism $\psi:P_1\to
P_2$ and an element $g\in G$ such that $h_{2}=gh_{1}g^{-1}$ and
$\Theta_2\circ\psi = g\cdot\Theta_1$.
 The leaves of $\pi^{\ast}{\fa}$ are the connected
components of $\Theta^{-1}(g),g\in G$ and we have a submersion
${\tilde \Theta}\colon M\rightarrow G/H$ such that $\sigma \circ
\Theta = \tilde \Theta \circ \pi$ where $\sigma\colon G\rightarrow
G/H$ is the quotient map.
\end{Proposition}

A proof of this proposition can be done (as in the classical way)
by constructing a suitable system of differential forms on $M$
which satisfy the same relations than the forms in a basis of the
Lie algebra $\mathcal L(G)$ of $G$ and then applying the classical
Darboux-Lie theorem and Ch.III, Cor. 2.4 and Prop. 2.7   in
\cite{Godbillon}.

\section{Construction of  a development}

\label{Section:development} In this section $\fa$ is a foliation
which is  $G$-i.u.t.a. on a connected manifold $M$. Our aim is to
give a self-contained proof of  a version of
Proposition~\ref{Proposition:Development} which is more adequate
to our approach and purposes. Indeed, we prove:

\begin{Proposition}
\label{PropositionDevelopment} Assume that $\fa$ is $G$-i.u.t.a.
Then $\fa$ has a $G$-transversal structure and a complementary
foliation. Moreover, given any leaf $L$ of $\fa$ we have a
development
of $\fa$ as follows: \\
\indent{\rm(i)} A Galoisian covering map $\pi\colon P=G\times
L\rightarrow
M$,\\
\indent{\rm(ii)} A homomorphism $h\colon \pi_{1}(M)\rightarrow G$
such
that $\pi^{-1}(e_G)\simeq Aut(\pi)$,\\
\indent{\rm(iii)} A submersion $\Theta\colon P\rightarrow G$ which
is a first integral for the pull-back foliation $\pi^{*}{\fa}$ and
is equivariant by $h$, that is,
$\Theta(\alpha.x)=h(\alpha)\Theta(x),\ \forall x\in M,\ \alpha \in
\pi_{1}(M)$.
\end{Proposition}

\begin{proof}[Proof of Proposition~\ref{PropositionDevelopment}]
Choose a leaf $L$ of  ${\cal F}$. The restriction $\Phi_L$ to
$G\times L$ of the action $\Phi \colon G\times M\to M$ has a
bijective tangent mapping at every point, which implies that it is
a local diffeomorphism, and that the image $\Phi(G\times L)$ is an
open subset of $M$. In $M$ we have the equivalence relation $x\sim
y$ if and only if there exists a $g\in G$ such that $\Phi(x,g)\in
L_y$, the leaf of $\fa$ through $y$, and the equivalence classes
are the open sets $\Phi(G\times L)$ where $L\in\fa$ is a leaf.
Since  $M$ is connected there only  one equivalence class, that is
$\Phi(G\times L)=M$. This shows that  the local diffeomorphism
$\Phi_L\colon G\times L \to M$ is surjective.

The assumption that the action of $G$ maps leaves of ${\cal F}$ to
leaves of ${\cal F}$ implies that for any $g\in G$ the following
conditions are equivalent: (a) There exists an $x\in L$ such
$\Phi(g,x)\in L$. (b) $g_M(L)=L$, if $g_M$ denotes the mapping
$x\mapsto\Phi(g,x)$. (The mapping $g\mapsto g_M$ is a homomorphism
from $G$ to the group of all diffeomorphisms of $M$.) Let $H=H_L$
denote the set of $g\in G$ which satisfy (a) or (b). Then $H$ is a
subgroup of $G$, and we have that $\Phi_L(g,x)=\Phi_L(g',x')$ if
and only if there exists an $h\in H$ such that $g'=gh^{-1}$ and
$x'=\Phi(h,x)$. The mapping $h\mapsto ((g,x)\mapsto
(gh^{-1},\Phi(h,x))$ defines an action of $H$ on $G\times L$,
which is free because the action on the first component is free.
The mapping $\Phi_L$ induces a bijective mapping $\Psi_L:G\times_H
L\to M$, which is the uniquely determined by the condition that
$\Phi_L = \Psi_L\circ p$ where $p$ denotes the canonical
projection from $G\times L$ onto the space $G\times_H L\to M$ of
$H$-orbits in $G\times L$.

The definition of $G\times_H L$ and $\Psi_L$ in the previous
paragraph was purely set-theoretic, let us  now discuss the
topological and smoothness aspects.

\begin{Claim}
\label{Claim:proper} If we provide $H$ with the discrete topology,
then the action of $H$ on $G\times L$ is a proper mapping, i.e.
the mapping $(h,(g,x))\mapsto ((gh^{-1},\Phi(h,x)), (g,x))$ is a
proper mapping from $H\times (G\times L)$ to $(G\times L)\times
(G\times L)$. \end{Claim}

\begin{proof}[Proof of Claim~\ref{Claim:proper}]
 Let us show that given an infinite sequence $h_j\in H$,
$g_j\in G$, $x_j\in L$, such that $(g_j h_j^{-1},\Phi(h_j,x_j))$
converges in $G\times L$ to $(g',x')$ and $(g_j,x_j)$ converges in
$G\times L$ to $(g,x)$, then a subsequence of the $(h_j,g_j,x_j)$
converges in $H\times G\times L$ to some element of $H\times
G\times L$. Here we use the leaf topology in $L$ (notice that this
is different from the $M$-topology in $L$ if $L$ is not a closed
subset of $M$). Because we use the discrete topology in $L$, this
amounts to the statement that $g_j\to g$ in $G$, $x_j\to x\in L$,
$g_j h_j^{-1}\to g'$ in $G$, $\Phi(h_j,x_j)\to x'$ in $L$ implies
that the $h_j$ have a constant subsequence.  From the fact that
the $g_j$ and the $g_j h_j^{-1}$ converge in $G$, it follows that
the $h_j = (g_j h_j^{-1})^{-1} g_j$ converge in $G$. There are
open neighborhoods $U$ and $V$ of $e_G$ and $x$ in $G$ and $L$,
respectively, such that $V$ is connected, and the restriction
$\psi$ of the $G$-action $\Phi$ to $U\times V$ is a diffeomorphism
from $U\times V$ onto an open neighborhood $W$ of $x$ in $M$, in
such a way that the integral manifolds of the restriction to $W$
of the vector subbundle $\mathbb {F}$ of $TM$, i.e.,  the ``local
leaves'', are equal to the sets of the form $\psi(\{u\}\times V)$,
where $u$ runs over $U$. That we have a diffeomorphism $\psi$
follows from the fact that the tangent mapping of $\psi$ at
$(e_G,x)$ is a bijective linear mapping from ${{\cal L}G}\times
T_xL$ onto $T_xM$, and the statement about the local leaves
follows from fact that the action maps leaves to leaves and
therefore the local action maps local leaves to local leaves. Let
$W_0$ be a closed neighborhood of $x$ in $M$ such that $W_0\subset
W$. Because the $h_j$ converge in $G$, and $x_j\to x$ in $M$,
there exist an integer $k$ such that $h_k^{-1} h_j\in U$, $x_j\in
V$ and $y_j:=\Phi(h_k^{-1} h_j,x_j)\in W_0$ whenever $j\geq k$.
Because for $j\to\infty$ the $y_j$ converge {\em in the leaf
topology} to an element of $W_0\subset W$, we conclude that there
exist a $k'$ such that all $y_j$ for $j\geq k'$, belong to the
same local leaf, which means that all $h_k^{-1} h_j$ are the same
for all $j\geq k'$, which in turn implies that the $h_j$ are the
same for all $j\geq k'$.
 This completes the proof of the claim.
 \end{proof}

Because the action of $H$ is proper and free, there is a unique
structure of smooth manifold on the orbit space $G\times_H L$ for
which the canonical projection $p: G\times L\to G\times_H L$ is a
principal $H$-bundle, in which $H$ is provided with the discrete
topology. In other words, the canonical projection $G\times L\to
G\times_H L$ is a Galois covering with group $H$. From the local
triviality we obtain that the mapping $\Psi_L:G\times_H L\to M$ is
smooth, because $\Phi_L = \Psi_L\circ p$, and $\Phi_L$ and $p$
were local diffeomorphisms, we obtain that $\Psi_L$ is a local
diffeomorphism, and because $\Psi_L$ is bijective, it follows that
$\Psi_L$ is a diffeomorphism from $G\times_H L$ onto $M$.

This completes the proof that $(P,h,\Theta)$, in which $P:=G\times
L$, $\pi := \Psi_L\circ p$ with covering group $H$, $h :=$ the
holonomy homomorphism from $\pi_1(M)$ onto the subgroup $H$ of $G$
and $\Theta :G\to G$ defined by $\Theta(g,x)=g$, is a development
of ${\cal F}$. The uniqueness of developments is straightforward.
\end{proof}

\section{Algebraic model and proof of
Theorem~\ref{Theoremfibrebundlenew}}
\label{Section:algebraicmodel}

Summarizing the discussion in the proof of
Proposition~\ref{PropositionDevelopment} we have: Let $\fa$ be a
foliation $G$-i.u.t.a. on $M$. Let $L$ be a leaf of the foliation
and let $H$ be the set of $g\in G$ such that $\Phi(g,l)\in L$ for
every (or equivalenty for some) $l\in L$. Then $H$ is a subgroup
of $G$, not necessarily closed, which we endow  the discrete
topology. Assuming that $M$ is connected, the restriction of
$\Phi$ to $G\times L$ is a covering map $\Phi :G\times L\to M$. We
have $\Phi(g,l)=\Phi(g',l')$ if and only if there is a uniquely
determined $h\in H$ such that $g'=gh^{-1}$ and $l'=\Phi(h,l)$. On
$G\times L$ we have the action of $H$ in which $h\in H$ sends
$(g,l)$ to $(gh^{-1},\Phi(h,l))$. This action is proper because
the action of $H$ on $L$ is proper (and discrete), and the action
is free because the right action of $H$ on $G$ is free. As a
consequence the action is proper and free, we have a smooth
quotient manifold $(G\times L)/H$, which is $G$-equivariantly
diffeomorphic to $M$, where the diffeomorphism from $(G\times
L)/H$ onto $M$ is induced by $\Phi$, and the $G$-equivariance is
with respect to the action of $G$ on $(G\times L)/H$ in which
$g'\in G$ sends the $H$-orbit of $(g,l)$ to the $H$-orbit of
$(g'g,l)$. If we write $p:G\times L\to (G\times L)/H$ for the
natural projection from $G\times L$ onto $(G\times L)/H$, then the
leaves of the foliation in $M$ are the sets $\Phi(p(\{g\}\times
L))$, $g\in G$. That is: \vglue.1in {\em The foliation $\fa$ in
$M$ corresponds to the foliation of the $p(\{g\}\times L)$, $g\in
G$, in $(G\times L)/H$.} \vglue.1in
 This is proves Theorem~\ref{Theorem:algebraicmodel}, i.e., the {\it algebraic model}
 $(G\times L)/H$ of the general foliation invariant under a
transverse Lie group action. In particular, all further analysis
can be done in $(G\times L)/H$, in which $H$ is a subgroup of $G$
acting from the right on $G$ and acting properly and discretely on
the manifold $L$.

Now we can prove several results.

\begin{Proposition}
\label{Proposition:closedleaf} The following statements are
equivalent.

\noindent{\rm(a)} The foliation has a closed leaf $L$.

\noindent{\rm(b)} $H$ is a closed discrete subgroup of $G$.

\noindent{\rm(c)} $H$ is a closed discrete subgroup of $G$ and the
projection $G\times L\to G$ onto the first factor exhibits
$(G\times L)/H$ as a fiber bundle over $G/H$ with fiber
diffeomorphic to $L$.
\end{Proposition}
\begin{proof}
 Suppose (a) holds and let $h_j$ be a sequence in $H$
which converges in $G$ to some $g\in G$. Let $x\in L$. Then
$\Phi(h_j,x)\to\Phi(g,x)$ in $M$ and therefore $\Phi(g,x)\in L$
because $L$ is closed in $M$. It follows that $g\in H$ proving
that (a) implies (b). That (b) implies (c) is a general fact
 about closed subgroups $H$ of a Lie group $G$, where
$H$ acts smoothly on a manifold $L$. Finally,  that (c) implies
(a) is obvious.
\end{proof}

Condition (c) means that the foliation in $M$ is a $G$-invariant
fibration. Thus we have:

\begin{proof}[Proof of Theorem~\ref{Theoremfibrebundlenew}]
The theorem follows from Proposition~\ref{Proposition:closedleaf}
and the above construction of the algebraic model.
\end{proof}

\begin{proof}[Proof of
Corollary~\ref{Corollary:finitehomotopy}] If $\pi_1(M,x)$ is
finite, then $H=h(\pi_1(M))$ is finite, hence a closed subgroup of
$G$, and the foliation is a fibration over the base space $G/H$.
If moreover $M$ is compact, then the fiber $L$ as well as the
group $G$ is compact. If in addition $G$ is simply-connected then
$G$ is semi-simple.
\end{proof}

\begin{proof}[Proof of Corollary~\ref{Corollary:codimensiontwo}]
In this case $G$ is a two-dimensional connected Lie group thus $G$
must be isomorphic to ${\mathbb R}^{2},\ S^{1}\times{\mathbb R},\
{\Aff}^{+}({\mathbb R})\simeq {\mathbb R}^{+}\times{\mathbb R}$ or
$S^{1}\times S^{1}$. The subgroup $H=h(\pi_{1}(M))$ must be
discrete.  If $H< {\mathbb R}^{2}$ then $H$ is isomorphic to the
trivial group, $0\times{\mathbb Z}$ or ${\mathbb Z}^{2}$, so
$G/H\simeq {\mathbb R}^{2},\ {\mathbb R}^{*}\times {\mathbb R}$ or
$T^{2}$. If $H<{S^{1}\times{\mathbb R}}$ then $H$ is isomorphic to
the trivial group, ${\mathbb Z}\times 0,\ 0\times{\mathbb Z}$ or
${\mathbb Z}^{2}$ so $G/H\simeq S^{1}\times{\mathbb R},\
{S^{1}\times{\mathbb R}\over{\mathbb Z}\times 0},\ S^{1}\times
S^{1},\ {S^{1}\times{\mathbb R}\over{\mathbb Z}\times{\mathbb
Z}}$. And if $H<{\Aff}^{+}({\mathbb R})$ then $H$ is isomorphic to
the trivial group, ${\mathbb N}\times 0,\ 0\times {\mathbb Z}$ or
${\mathbb N}\times{\mathbb Z}$. Since $M$ is compact, $G/H$ is
compact therefore we have $G/H\simeq T^{2}$.
\end{proof}

\section{Foliations invariant under a transverse local action}

\label{Section:localaction} Proposition~\ref{Proposition:Lieforms}
states that the existence  of a $G$-transversal structure for a
foliation $\fa$ is equivalent to the existence of a suitable
system of differential forms $\{\Omega_j\}_{j=1}^q$ satisfying the
same structure equations of a given basis of the Lie algebra
$\mathcal L(G)$. Let us now prove this and give an interpretation
of the invariance of $\fa$ under a $G$-transversal action  in a
way that motivates a generalization. Thus, on what follows we
assume that $\fa$ is a foliation on $M$ which is $G$-i.u.t.a.  For
any $X\in {{\cal L}G}$, let $X_M$ denote vector field on $M$ which
defines the infinitesimal action of $X$ on $M$. The assumption
that the $G$-orbits have the same dimension as $G$ is equivalent
to the condition that the action is locally free, which in turn is
equivalent to the condition that for each $x\in M$ the mapping
$X\mapsto X_{M,x}$ is injective from ${{\cal L}G}$ to $T_xM$.
Denote the image space by ${{\cal L}G}_{M,x}$, this can be viewed
as the tangent space at $x$ to the orbit through $x$. The
transversality condition means that, for any $x\in M$, $T_xM$ is
equal to the direct sum of $\mathbb {F}_x$ and ${{\cal
L}G}_{M,x}$. Therefore,  there is a unique ${{\cal L}G}$-valued
one-form $\Theta$ on $M$, such that $\Theta=0$ on $\mathbb {F}$
and $\Theta_x(X_{M,x})=X$ for every $X\in {{\cal L}G}$. This
${{\cal L}G}$-valued one-form $\Theta$ on $M$ is called the {\it
connection form} of the infinitesimal connection $\mathbb {F}$ for
the infinitesimal action of ${{\cal L}G}$ on $M$. The form
$\Theta$ is automatically smooth. The fact that the infinitesimal
connection is flat, meaning that $\mathbb {F}$ is integrable, is
equivalent to the condition that $(d\Theta)(u,v) =
-[\Theta(u),\Theta(v)]$ for every pair of vector fields $u$, $v$
on $M$. The one-forms $\Omega_j$ appearing in
Proposition~\ref{Proposition:Development} are exactly the
components of the connection form. Thus we can prove
Proposition~\ref{Proposition:Lieforms} by a repetition of the
proof that integrability of $\mathbb {F}$ implies that
$(d\Theta)(u,v)= -[\Theta(u),\Theta(v)]$, or equivalently $d
\Omega_k=\sum_{i<j} c^k_{ij}\Omega_i\wedge\Omega_j$. Now, a
construction of the mappings $f_i$ and $\gamma_{ij}$ in the
definition of Lie transverse structure
(Definition~\ref{Definition:Liestructure}) is not clear at a first
moment. This is quite obvious  once we observe that this is
equivalent to a construction of a local action of $G$ on the leaf
space and that  the local action of $G$  maps local leaves to
local leaves.

 Motivated by this we observe that a weaker assumption
 than the invariance under a Lie group transverse action
 is the following:

\begin{Definition}
\label{Definition:invariantlocalaction} {\rm We  say that a
foliation $\fa$ in $M$ is {\it invariant under a transverse local
action} of a Lie group $G$ ($G$-i.u.t.l.a.)  if there is a locally
free local action of $G$ on $M$, the tangent mappings of which
leave $\mathbb {F}$ invariant, and such that the ${{\cal
L}G}_{M,x}$ are complementary subspaces to the $\mathbb {F}_x$ in
$T_xM$.}
\end{Definition}

The relation between the two notions in
Definitions~\ref{Definition:Liestructure} and
\ref{Definition:invariantlocalaction} is given below:

\begin{Proposition}
Given $\fa$ on $M$ the following conditions are equivalent:

\noindent{\rm (i)} $\fa$ is $G$-i.u.t.l.a.

\noindent{\rm (ii)} $\fa$ has a transversal $G$-structure and a
complementary foliation.
\end{Proposition}

\begin{proof} Assume that $\fa$ is $G$-i.u.t.l.a.
Because the ${{\cal L}G}_{M,x}$ are the tangent spaces to
the local $G$-orbits, they define an integrable vector subbundle
${{\cal L}G}_M$ of $TM$, which is complementary to $\mathbb {F}$.
As we have already observed above  this weaker condition already
implies that there is a $G$-transverse structure. Conversely, if
$K$ is an integrable vector subbundle of $TM$ (i.e. defining a
foliation in $M$) which is complementary to $\mathbb {F}$, and we
have a transversal $G$-structure to $\mathbb {F}$, defined by a
${{\cal L}G}$-valued one-form $\omega$ as in 1), then for each
$X\in{{\cal L}G}$ and $x\in M$ there is a unique $X_{M,x}\in K_x$
such that $\omega_x(X_{M,x})=X$. This defines a smooth vector
field $X_M$ on $M$ and the equation $(d\omega)(u,v)=
-[\omega(u),\omega(v)]$ in combination with the integrability of
$K$ implies that $x\mapsto X_M$ is a homomorphism of Lie algebras
from ${{\cal L}G}$ to the Lie algebra of smooth vector fields in
$M$. In other words, in this way we obtain an infinitesimal, and
hence a local action of $G$, which also maps local leaves to local
leaves. It is the unique infinitesimal action of ${{\cal L}G}$ for
which $\omega$ is equal to the connection form and $K$ is tangent
to the orbits.
\end{proof}

{\small \begin{Remark} {\rm The assumption of having a transverse
$G$-action which maps leaves of ${\cal F}$ to leaves of ${\cal F}$
is equivalent to the weaker assumption that $\fa$ is
$G$-i.u.t.l.a., but with the additional assumption that the local
action of $G$ on $M$ can be extended to a global one on $M$. (If
such an extension exists, it is unique.) If $M$ is compact, an
extension to a global action always exists. Therefore, if $M$ is
compact, then $\fa$ is $G$-i.u.t.a. iff $\fa$ is $G$-i.u.t.l.a.,
i.e., the weaker assumption in
Definition~\ref{Definition:invariantlocalaction}  is equivalent to
the fact that $\fa$ is $G$-i.u.t.a. .

 Clarifying  how much
the assumption ``$\fa$ is $G$-i.u.t.a."
 is stronger than ``$\fa$ is $G$-i.u.t.l.a."  might help in
understanding which consequences are typical consequences of the
first and not only of the existence of a $G$-transverse structure.
In a rough manner,  a $G$-transversal structure for the foliation
${\cal F}$ is something like a locally free local action of $G$ on
the leaf space $M/{\cal F}$, in which the latter has to be treated
as a sort of non-Hausdorff manifold.}
\end{Remark}

}

\section{Solvable groups}
\label{Section:solvablegroups}

Let us consider the Lie group of affine maps of the real line
$\Aff({\mathbb R})=\{ x\mapsto a.x+b;\ a\in{\mathbb R}^{*}\mbox{
and } b\in{\mathbb R}\}\simeq {\mathbb R}^{*}\times{\mathbb R}$.
Let $\fa$ be a foliation of codimension two on $M$ invariant under
a transverse action of $\Aff(\mathbb R)$. We assume that $\fa$ is
transversely oriented so that indeed $\fa$ is ${\Aff}^{+}(\mathbb
R)$-i.u.t.a., where ${\Aff}^{+}(\mathbb R)$ is the subgroup of
orientation preserving affine maps of the real line. Notice that
as a manifold we have ${\Aff}^{+}(\mathbb R)=(0,+\infty) \times
\mathbb R$ so that it is simply-connected, also it is solvable as
a group. According to Proposition \ref{Proposition:Lieforms} there
is an integrable system of two one-forms $\omega,\eta$ defining
$\fa$ on $M$ such that $d\omega=\omega \wedge \eta, \, d\eta=0$.
If $H_{1}(M,{\mathbb R})=0$ then $\eta=dh$, for some
differentiable function $h\colon M\rightarrow {\mathbb R}$, and we
define $f=e^{h}\colon M\rightarrow{\mathbb R}^{*}$, thus
$\eta={1\over f}df$. A straightforward computation then shows that
$d(f{\omega})=0$ and therefore   $\omega =f^{-1}dg$ for some
differentiable function $g\colon M\rightarrow {\mathbb R}$. So
there exists a fibration ${F}=(f,g)\colon M\to \Aff({\mathbb R})$
whose fibers are the leaves of $\fa$. Assume now that $M$ is
compact. In this case according to Tischler's Theorem
\cite{Tischler}, since $\eta$ is a nonsingular closed one-form in
$M$, there exists a fibration $f\colon M\to S^{1}$. We have
proved:

\begin{Proposition}
\label{Proposition:affine} Let $M$ be a connected manifold with a
foliation $\fa$ which is  $\Aff({\mathbb R})$-i.u.t.a. Then $\fa$
has an $\Aff({\mathbb R})$-transverse structure and we have:
\begin{itemize}
\item[{\rm (i)}] If $H^{1}(M,{\mathbb R})=0$ then $\fa$ is given
by a submersion $F\colon M\to \Aff({\mathbb R})$. In particular,
in this case $M$ is not compact. \item[{\rm (ii)}] If $M$ is
compact then it admits a fibration $f\colon M\to S^{1}$.
\end{itemize}
\end{Proposition}

This proposition is a particular case of
Theorem~\ref{Theorem:solvable} of which  proof is given below:

\begin{proof}[Proof of Theorem~\ref{Theorem:solvable}]
 Let $D{{\cal L}G}$ denote the {\it derived Lie algebra} of $\mathcal L(G)$, i.e.,
 the  linear
subspace of ${{\cal L}G}$ which is generated by all $[X,Y]$ such
that $X, Y\in {{\cal L}G}$. It is  a known that $D\mathcal L (G)$
 is a normal Lie subalgebra of ${{\cal L}G}$ (\cite{JJ}). Let $\xi\in
(D{{\cal L}G})^0\simeq ({{\cal L}G}/D{{\cal L}G})^*$ be a linear
form on ${{\cal L}G}$ which is equal to zero on $D{{\cal L}G}$. If
$\Theta$ denotes the connection form introduced in
Section~\ref{Section:localaction} then $\omega := \omega_{\xi} :=
\xi\circ\Theta$ is a closed one-form on $M$, hence $\omega = df$
for a smooth real-valued function $f=f_{\xi}$ on $M$, which we can
let depend linearly on $\xi$. It follows that $x\mapsto
(\xi\mapsto f_{\xi}(x))$ defines a smooth mapping $f$ from $M$ to
$(({{\cal L}G}/D{{\cal L}G})^*)^* = {{\cal L}G}/D{{\cal L}G}$. It
is a submersion, the leaves of ${\cal F}$ are contained in the
fibers of $\mathbb {F}$, and also the fibers of $\mathbb {F}$ are
invariant under the action of the group $DG$ generated by $D{{\cal
L}G}$. According to Theorem 3.18.1 in \cite{V} if $G$ is simply
connected, the analytic subgroup $DG$ of $G$ defined by $D{{\cal
L}G}$ is a closed normal subgroup of $G$, therefore $G/DG$ is an
Abelian and simply connected Lie group.

The Abelian group $G/DG$ is isomorphic with ${{\cal L}G}/D{{\cal
L}G}$ and acts on it by translations. Using lifts by elements of
$G$ acting on $M$, we obtain that $f$ is surjective and defines a
topologically trivial fibration. In particular the first
cohomology group $H^1$ of each fiber of $f$ is equal to zero as
well. It is also clear that an element $g$ of $G$ leaves a fiber
of $f$ fixed if and only $g\in DG$. Because any $h\in H$ leaves
$L$ fixed, and therefore also leaves the $f$-fiber containing $L$
fixed, we conclude that $h\in DG$, i.e. we obtain that $H\subset
DG$.

This means that in the fibers of $f$, we have the same situation
again, with $H\subset DG$, connected $f$-fibers and the $H^1$ of
the $f$-fibers equal to zero. Since the Lie algebra ${{\cal L}G}$
is {\em solvable} (because $G$ is solvable), then the repeated
derived Lie algebras $D^i{{\cal L}G}$ terminate at zero (cf. p.
201 in \cite{V}), and we arrive at the conclusion that the group
$H$ is trivial. In view of Theorem~\ref{Theoremfibrebundlenew}
this means that $M$ is isomorphic to $G\times L$ and the foliation
is defined by the projection onto the first factor. Using Theorem
3.18.11 in \cite{V} we conclude that (since $G$ is simply
connected and solvable) $G$ it is diffeomorphic to an Euclidean
space. In the above we have  used the fact that any fiber bundle
over a contractible space is trivial. This can be found in
Corollary 11.6 in \cite{Steenrod}.
\end{proof}

We have the following generalization of
Theorem~\ref{Theorem:solvable}:

\begin{Theorem}
Let ${\cal F}$ be a smooth foliation in a connected smooth
manifold $M$. Assume that ${\cal F}$ is invariant under a
transverse  action of a Lie group $G$ and assume that $H^1(M,
\mathbb R)=0$. Let $p$ be the smallest nonnegative integer such
that $D^{p+1}{{\cal L}G}=D^p{{\cal L}G}$ and let $D^pG$ denote the
analytic Lie subgroup of $G$ with Lie algebra equal to $D^p{{\cal
L}G}$. Then $H\subset D^pG$.
\end{Theorem}

\begin{proof}  By passing to the universal covering of $G$,
we may assume that $G$ is simply connected. Theorem 3.18.12 in
\cite{V} states that in this case every analytic subgroup of $G$
is closed and simply connected. Then the result follows from the
proof of Theorem~\ref{Theorem:solvable}.
\end{proof}


\begin{proof}[Proof of Corollary~\ref{Corollary:d-torus}]
 We shall use the same notation of the proof of Theorem~\ref{Theorem:solvable}.
 Let $M$ be a compact smooth
manifold and ${\cal F}$ a foliation in $M$ which admits a
$G$-transversal structure. Then the mapping $\xi\mapsto
[\xi\circ\Theta]$ from $({{\cal L}G}/D{{\cal L}G})^*$ to $H^1_{de
Rham}(M)$ is injective. Moreover, if $d:=\dim H^1_{de Rham}(M)\leq
\dim {{\cal L}G} - \dim D{{\cal L}G}$, then we have equality here
and $M$ fibers over a $d$-dimensional torus in such a way that the
leaves of ${\cal F}$ are contained in the fibers of this
fibration. If ${{\cal L}G}$ is abelian and $\dim H^1_{de
Rham}(M)\leq \dim {{\cal L}G}$, then ${\cal F}$ is a fibration
over a torus.
\end{proof}

\section{Holomorphic foliations}
\label{section:Holomorphic}

In this section we study holomorphic foliations which are
invariant under transverse actions of complex Lie groups.

\begin{Theorem}\label{Theorem:locallyholomorphictrivialfibration}
Let $M$ be a  connected complex manifold and $\fa$ a holomorphic
foliation invariant under a holomorphic transverse action of a
complex Lie group $G$ of dimension $\dim G=\codim \fa$. The
following conditions are equivalent:

\begin{itemize}
\item[{\rm (a)}] ${\cal F}$ has a leaf $L$ which is closed in $M$.

\item[{\rm (b)}] $H(L)$ is a discrete  subgroup of $G$.

\item[{\rm (c)}] The projection $\pi:G\times L\to G$ onto the
first factor induces a {\underline {holomorphic}}  fibration
$M\simeq G\times_H(L) L\to G/H(L)$, of which the fibers are the
leaves of ${\cal F}$.
\end{itemize}
\end{Theorem}

\begin{Remark}{\rm In the above statement the fibration
$M\simeq G\times_H(L) L\to G/H(L)$ is a {\it holomorphic}
fibration, in the sense that it the local trivializations are
biholomorphic maps.  According to Ehresmann's Theorem
\cite{Godbillon}, any proper $C^r, r \geq 2$ submersion defines a
$C^r$-locally trivial fiber bundle. This is not true for proper
holomorphic submersions. Indeed, the analytic type of the fiber
may vary. On the other hand, Grauert-Fischer's theorem
\cite{Barth} asserts that this is the only obstruction: A proper
holomorphic  submersion  is a holomorphic fibration (i.e., a
locally trivial holomorphic fiber bundle) if and only if the
fibers  are holomorphically equivalent. Thus
Theorem~\ref{Theorem:locallyholomorphictrivialfibration} follows
from Theorem~\ref{Theoremfibrebundlenew} and Grauert-Fischer'
theorem. Nevertheless, we give a ``simpler" self-contained proof
in Section~\ref{section:Holomorphic}.}
\end{Remark}

\begin{proof} [Proof of
Theorem~\ref{Theorem:locallyholomorphictrivialfibration}]
 We already know (Theorem~\ref{Theoremfibrebundlenew}) that $M$ is a
$C^\infty$ fibre bundle over the homogenous space $G/H$. Since the
fibers are holomorphically equivalent (by biholomorphisms
$\Phi_{g}\colon M\rightarrow M$), Grauert-Fischer's Theorem
\cite{Barth} states that the submersion $M\rightarrow G/H$ is a
locally holomorphically trivial fibration. Nonetheless, we can
give a more self-contained proof as follows:  We know that $H$ is
a closed and zero-dimensional subgroup of $G$ and that $\pi :M\to
G/H$ is a $G$-equivariant holomorphic mapping. Let $g\in G$ and
write $L=\pi^{-1}(\{gH\})$ for the fiber over the point $gH$ in
$G/H$. Then there exists an open neighborhood $U$ of the origin
 in ${{\cal L}G}$ such that $X\mapsto \exp(X)gH$ is a
holomorphic diffeomorphism from $U$ onto an open neighborhood $V$
of $gH$ in $G/H$. The mapping $(X,x)\mapsto\Phi(\exp(X),x)$ is a
holomorphic diffeomorphism from $U\times L$ onto $\pi^{-1}(V)$,
and it yields the desired holomorphic trivialization.
\end{proof}

\begin{Proposition}
\label{Proposition:KahlerLiestructure} Let $M$ be a compact
K\"ahler manifold, $G$ a complex simply connected Lie group, and
${\cal F}$ a holomorphic foliation of codimension $q$ with a
{\rm(}holomorphic{\rm)} $G$-transverse structure. Then $G\simeq
\C^q$. If moreover $\dim H^1(M, \mathbb R)\leq 2q$, then $\dim
H^1(M, \mathbb R)= 2q$ and the foliation ${\cal F}$ is a fibration
over a real $2q$-dimensional torus.
\end{Proposition}

\begin{proof}
According to \cite{Gr-Ha} p. 110 any holomorphic $q$-form on a
compact K\"ahler manifold is closed. Applying this to the
connection form  (which is holomorphic), we conclude that this is
closed, which in turn implies that ${{\cal L}G}$ and therefore $G$
is abelian.
\end{proof}

A natural holomorphic version of
Proposition~\ref{PropositionDevelopment}  implies the following:

\begin{Proposition}\label{Proposition:complexcovering} Let $M$ be
a compact K\"ahler manifold with a holomorphic codimension $q$
foliation $\fa$ invariant under a Lie group transverse action of
$G$. Then the universal covering of $G$ is isomorphic to
$(\C^{q},+)$. If moreover $\fa$ has a compact leaf then we have
$G=\C^{q}/H$ for some closed subgroup $H < \C^{q}$.
\end{Proposition}

{\small\begin{Remark} {\rm Since an algebraic manifold is  always
K\"ahler, Proposition~\ref{Proposition:complexcovering} is valid
for any projective manifold.}
\end{Remark}}

\subsubsection*{Codimension one algebraic foliations}  By definition such an
{\it algebraic foliation} $\fa_0$  on $\C^n$  is   given by a
polynomial one-form $\Omega= \sum\limits_{j=1}^n P_{j} dz_{j}$,
where the $P_j$ are polynomials in the affine variables
$(z_1,...,z_n)\in \C^n$, satisfying the integrability condition
$\Omega \wedge d \Omega=0$. Such a foliation admits a unique
extension to a holomorphic foliation with singularities $\fa$ on
$\C\mathbb P^{n}$. Conversely, any foliation of $\C\mathbb P^{n}$
is obtained this way. Assume now that $\fa_0$ is $\C$-i.u.t.a.,
i.e., invariant by a holomorphic flow in $\C^n$. The foliation is
then given by a closed holomorphic one-form $\omega$ on $\mathbb
C^n$. Thus we have $\omega= dF$ for an entire function ${F}$ on
$\mathbb C^n$.

\begin{Claim}\label{Claim:meromorphicextension}
If the hyperplane $\C\mathbb P^{n-1}_\infty=\C\mathbb
P^{n}\setminus \C^n$ is not $\fa$-invariant then ${F}$ is a
polynomial first integral on $\C\mathbb P^{n}$.
\end{Claim}
\begin{proof}
 Fixed a generic point $q \in \C\mathbb P^{n-1}_\infty$  we may
consider a ``flow box" (i.e., a distinguished neighborhood for
$\fa$) \, $V$ containing $q$ with coordinates
$(z_{1},...,z_{n})\in V $, such that  $\C\mathbb
P^{n-1}_\infty\cap V= \{z_1=0\}$ and $\fa\big |_V$ is given by
$dz_n=0$. Let $V ^* = V \setminus (V \cap \C P(n-1)_\infty)= V
\setminus \{z_1=0 \}$. In $V^*$ we have $\omega\wedge dz_n\equiv
0$, i.e., $dF\wedge dz_n\equiv 0$. Therefore $
{F}\big|_{V^*}=F(z_n)$ is depends only on the variable $z_n$. On
the other hand, $ {F}\big|_{V^*}$ is holomorphic. Therefore ${F}$
extends meromorphically to $V$. Then Hartogs' theorem
\cite{Gunning1} implies that $ {F}$ is meromorphic on $\C\mathbb
P^{n}$. Liouville's theorem \cite{Gunning2} then shows that ${F}$
is a rational function and since it is holomorphic on $\C ^n$ we
conclude that $ {F}$ is a polynomial on $\C^n$.\end{proof}

\begin{Proposition}
Let $\fa_0$ be an algebraic codimension one foliation on $\C^n, n
\geq 2$. Suppose that $\fa$ is $\C$-i.u.t.a. Then the hyperplane
at infinity $\C P(n-1)_\infty$ is $\fa$-invariant.
\end{Proposition}
\begin{proof} If $\C \mathbb P^{n-1}_\infty$ is not invariant then by the above
claim $\fa$ has a polynomial first integral $F$  on $\C^n$.
However, as a general fact for meromorphic first integrals, the
polar set $\{F=\infty\}$ and the zero set $\{F=0\}$ are invariant.
Since  the polar set is $\C \mathbb P^{n-1}_\infty$ the
proposition follows.
\end{proof}

Using techniques introduced by R. Mol in \cite{Mol} one may be
able to go further in the classification of $\fa$ in this case.

\subsubsection*{Codimension-$q$ foliations on complex projective spaces}

 Let $\fa$ be a codimension one holomorphic foliation (necessarily with singularities)
 on $\C\mathbb P^{n}$. Suppose that we have an automorphism
$\Phi\colon\C^{n}\rightarrow \C^{n}$ such that
$\Phi^{*}{\fa}_{0}\equiv {\fa}_{0}$ where $\fa_0$ is the
restriction of $\fa$ to $\C^n$.

\begin{Lemma}\label{Lemma:autoextends}If $\C\mathbb P^{n-1}_\infty$ is not $\fa$-invariant,
then $\Phi$ extends meromorphically to $\C\mathbb P^{n}$ and
therefore $\Phi$ is algebraic.
\end{Lemma}
\begin{proof} Let $L\in{\fa}$ be a generic leaf, so that the holonomy
group $\Hol(L)$ is trivial. We fix a ``flow box" $U$ containing a
point $q_{0}\in L\cap\C\mathbb P^{n-1}_\infty$ and a transverse
disk $\Sigma$ centered at a point  $q\in (L\cap U)\setminus
\C\mathbb P^{n-1}_\infty$ near $q_{0}$. Next, for $L_{1}=\Phi(L)$,
we fix a ``flow box" $V$ containing a  point $p_{0}\in L_{1}\cap
\C\mathbb P^{n-1}_\infty$ and a transverse disk $\widetilde
\Sigma$ for $p\in (L_{1}\cap V)\setminus \C\mathbb P^{n-1}_\infty$
near $p_{0}$. Finally, we consider $\Sigma_{1}=\Phi(\Sigma),\
q_{1}=\Phi(q_{0})\in L_{1}$ and a path $\alpha$ from $q_{1}$ to
$p$ in the  leaf $L_{1}$: for a point $z\in \Sigma$ there exists a
unique point ${\widehat z}\in \C\mathbb P^{n-1}_\infty$ such that
$z$ and ${\widehat z}$ are in the same plaque of $\fa$ on $U$,
thus we may define a map $f\colon \Sigma \rightarrow
U\cap\C\mathbb P^{n-1}_\infty$ by $f(z)={\widehat z}$. In the same
way, we define $g\colon \Sigma_{1}\to V\cap\C\mathbb
P^{n-1}_\infty$. By the holonomy map $h_{\alpha}$ induced by
$\alpha$, we define $\psi\colon \C\mathbb
P^{n-1}_\infty\rightarrow \C\mathbb P^{n-1}_\infty$ by
$\psi({\widehat z})=g\circ h_{\alpha}\circ \Phi \circ
f^{-1}({\widehat z})$. Since $\Hol(L_{1})=0$, $h_{\alpha}$ is
unique and $\psi$ is well defined. The restrictions ${\fa}|_{U}$
and ${\fa}|_{V}$ are trivial so, for any leaf with trivial
holonomy, $\psi$ is an extension of $\Phi$ to $\C\mathbb P^{n}$.
By Hartogs' Extension Theorem, $\Phi$ extends to a map $\Phi\colon
\C\mathbb P^{n}\rightarrow\C\mathbb P^{n}$. Since the inverse
$\Phi^{-1}$ also extends to ${\C\mathbb P^{n}}$ we conclude that
$\Phi$ is an automorphism of $\C\mathbb P^{n}$. \end{proof}

These very same ideas give:

\begin{Theorem}\label{Theorem:extension} Let $\fa$ be a
codimension $q$ singular holomorphic foliation on $\C\mathbb
P^{n}$, $\Gamma \subset \C\mathbb P^{n}$ an algebraic irreducible
hypersurface which is not $\fa$-invariant and $\Phi\colon
{\C\mathbb P^{n}\setminus \Gamma} \rightarrow{\C\mathbb
P^{n}\setminus \Gamma}$ a holomorphic diffeomorphism preserving
the foliation. Then $\Phi$ is the restriction of a birational map
of $\C\mathbb P^{n}$.
\end{Theorem}

\begin{proof}[Proof of Theorem~\ref{Theorem:extension}]
 Let $\fa$ be a codimension $q$ holomorphic foliation on $\C\mathbb P^{n}$  and
$\Phi\colon \C\mathbb P^{n}\setminus \Gamma \rightarrow\C\mathbb
P^{n}\setminus \Gamma$ a holomorphic diffeomorphism preserving
$\fa$. Let $L$ be a leaf of $\fa$ transverse to $\Gamma$.  We have
$\dim (L\cap\Gamma)=(n-q+n-1)-n=n-q-1$, so we consider a
$(n-q-1)-$disk $\Sigma$ transverse to $L$ and we obtain, as in the
proof of Lemma~\ref{Lemma:autoextends}, that the automorphism
extends to a neighborhood of a point $q\in \Sigma \cap \Gamma$ in
$\C\mathbb P^{n}$ and therefore to $\C\mathbb P^{n}$.
\end{proof}

From the above result we promptly obtain
Theorem~\ref{Theorem:extensionaction}.
Theorem~\ref{Theorem:extensionaction} and
Proposition~\ref{Proposition:complexcovering} then give the
description of the foliation in
Theorem~\ref{Theorem:extensionaction}.

{\small \begin{Remark}[Singular foliations] {\rm  We consider a
singular codimension $q\ge 1$ holomorphic foliation ${\fa}$ on
complex manifold $M$. We shall always assume that the singular set
$\sing({\fa})$ of $\fa$ has codimension $\geq 2$. Denote by $\fa
^\prime$ the underlying nonsingular foliation on $M\setminus
\sing(\fa)$. Let $\Phi\colon G\times M\rightarrow M$ be a
holomorphic action of a complex Lie group $G$. We say that $\fa$
{\it is invariant under the  transverse action $\Phi$} if (1)
$\Phi_g(\sing({\fa}))=\sing({\fa})$, for all $g \in G$, and
$\fa^\prime$ is $G$-i.u.t.a. with respect to $\Phi$.  When $q=1$,
we have a flow (say given by a complete holomorphic vector field
$X$ on $M$) under which  $\fa$ and $\sing({\fa})$ are invariant.
There exists a holomorphic closed one-form $\omega \in
\Lambda^{1}(M\setminus \sing({\fa}))$ which defines the foliation.
Since $\Cod (\sing({\fa}))\geq 2$ Hartogs' Extension Theorem
\cite{Gunning1} implies that the one-form $\omega$ extends
holomorphically to $M$. We conclude that $\fa$ {\it is
nonsingular} because we cannot have $\omega \cdot X\equiv 1$ on
$M$ if $\omega$ is holomorphic and $X$ has singularities. This
suggests  that  the interesting case occurs when the foliation
admits a Lie group transverse action in the complement of a
codimension one invariant analytic subset $\Lambda$ such that
$\sing(\fa)\subset \Lambda$. This is the case of linear foliations
on complex projective spaces (\cite{Scardua}).}
\end{Remark}}

\section{Complements}

\subsection{The Realization problem}

In \cite{Meigniez} the author discusses, mainly for solvable Lie
groups,  the ``{\em Realization problem}'', which is the following
question of Haefliger:

\begin{Question} Which subgroups $H$ of a given Lie group $G$ can occur as
$h(\pi_1(M))$ for a development $(P,h,\Theta)$ of a $G$-transverse
structure on a compact manifold $M$?
\end{Question}

Under the additional hypothesis that the  ``$G$-transverse
structure'' is followed by the fact that the  ``the foliation is
invariant under a transverse action of $G$ on $M$'', then the
realization problem asks for the subgroups $H$ of $G$ for which
there exists a smooth connected manifold $L$ such that $H$ acts
smoothly on $L$, the right-left-action of $H$ on $G\times L$ is
proper (here $H$ is provided with the discrete topology) and the
quotient $G\times_HL$ is compact. This seems to imply that $H$ is
finitely generated. Also, for closed subgroups $H$ of $G$ the
answer is that $H$ is a discrete subgroup of $G$, $G/H$ is
compact, and one can take any smooth action of $H$ (for instance,
the trivial action) on any compact connected manifold $L$.
Nevertheless, in Theory of Foliations  one is mostly interested in
foliations which are not fibrations, which corresponds to the case
that $H$ is not closed in $G$. Such cases can occur, as shown by
the example of orbit foliations of non-closed subgroups of tori.
 It may be  helpful   to restrict
to solvable Lie groups $G$.

\subsection{More on the algebraic model}
Let us say a few more words about  the algebraic model $G\times_H
L$ for a foliated manifold with transverse Lie group action. The
main motivation comes from the book \cite{JJ}, more precisely from
its  Section 2.4 where the authors introduce the associated fiber
bundle $X\times_H Y$, in which $X$ and $Y$ are manifolds, $H$ is a
Lie group acting on $X$ by $(h,x)\mapsto x\cdot h^{-1}$, on $Y$ by
$(h,y)\mapsto h\cdot y$, and it is assumed that the action of $H$
on $X$ is proper and free. This then implies that the action
$(h,(x,y))\mapsto (x\cdot h^{-1}, h\cdot y)$ of $H$ on $X\times Y$
is proper and free, which then makes that the orbit space
$X\times_HY$ is a smooth manifold, and that the projection
$X\times Y\to X\times_H Y$ is a principal $H$-bundle. Moreover,
the projection onto the first factor induces a fibration
$X\times_H Y\to X/H$, with fibers isomorphic to $Y$. In the
statements that the orbit spaces = quotients are smooth manifolds,
the properness assumption is quite essential.

If $X=G$ is equal to a Lie group and $H$ is a Lie subgroup of $G$,
then the right action of $H$ is proper if and only if $H$ is
closed in $G$. In this case there is a unique action of $G$ on
$G\times_H Y$ such that the projection $G\times Y\to G\times_H Y$
intertwines the left action of $G$ on $G\times Y$ (defined by
means of the left action of $G$ on the first factor) with the
action of $G$ on $G\times_H Y$. Furthermore the projection
$G\times _HY\to G/H$ intertwines the action of $G$ on $G\times_H
Y$ with the transitive left action of $G$ on $G/H$. This leads to
a well-understood model of a $G$-homogeneous bundle over a
$G$-homogeneous base manifold.

The point of this construction in the aforementioned book was that
for any proper Lie group action every orbit has an invariant open
neighborhood which, as a manifold with Lie group action, is
isomorphic to some $G\times_H Y$, where $H$ is the stabilizer
subgroup of a point in the given orbit and $Y$ is a so-called
slice for the $G$-action. This is the Tube theorem 2.4.1 in the
book.

It is then clear that set-theoretically the $G$-space $M$ is of
the form $G\times_H L$, in which $L$ is a leaf of your foliation
and $H$ is the subgroup of $G$ which maps $L$ to itself. Moreover,
 in order that $G\times_H L$ inherits the
structure of a smooth manifold from its construction as the space
of $H$-orbits (which is equivalent to saying that the mapping
$G\times L\to M$, defined by the restriction to $G\times L$ of the
$G$-action in $M$, is a principal $H$-bundle), it is sufficient to
assume that the $H$-action on $G\times L$ is proper. Finally, in
our situation  this properness follows from your assumptions (in
the more general situation that the $G$-action maps leaves to
leaves and is transversal, but $X_M$ tangent to a leaf for nonzero
elements $X$ of the Lie algebra, then one would need quite
technical additional assumptions on the action in order to obtain
that the action of $H$ on $G\times L$ is proper). Because $H$ is
discrete in our case, the principal $H$-fibration $G\times L\to M$
is a Galois covering (actually Galois coverings are nothing else
than principal fiber bundles with discrete group actions).

From the background sketched above
Theorem~\ref{Theoremfibrebundlenew}  is clear, i.e.,  $L$ is
closed in $M$ if and only if $H$ is closed in $G$ if and only if
$M$ is a $G$-homogeneous fiber bundle over the homogeneous space
$G/H$, where the fibers are equal to the leaves of the foliation.

\bibliographystyle{amsalpha}

\begin{thebibliography}{31}
\frenchspacing

\bibitem{Barth} W. Barth, C. Peters, A. Van de Ven: Compact
Complex Surfaces, Springer-Verlag, Berlin, 1984.


\bibitem{Godbillon} C. Godbillon; Feuilletages: {\'E}tudes
Geom{\'e}triques, Birkh{\"a}user, Berlin, 1991.


\bibitem{Gr-Ha} P. Griffiths; J. Harris: Principles of algebraic geometry.
Pure and Applied Mathematics. Wiley-Interscience [John Wiley \&
Sons], New York, 1978.

\bibitem{Gunning1} R.C. Gunning: Introduction to holomorphic
functions of several variables, vol. I, Function Theory, Wadsworth
and Brooks/Cole, Pacific Grove, CA, 1990.

\bibitem{Gunning2} R.C. Gunning: Introduction to holomorphic
functions of several variables. Vol. II. Local theory.  Wadsworth
and Brooks/Cole, Monterey, CA, 1990.

\bibitem{Meigniez}  G. Meigniez: {\it Holonomy
groups of solvable Lie foliations}. Integrable systems and
foliations/Feuilletages et syst\`{e}mes int\'{e}grables
(Montpellier, 1995), 107--146, Progr. Math., 145, Birkh\"{a}user
Boston, Boston, MA, 1997.



\bibitem{Mol} R. S. Mol: {\it Meromorphic first integrals: some extension results};
 Tohoku Math. J. (2) 54 (2002), no. 1, 85--104.


\bibitem{JJ} J.J. Duistermaat, J.A.C. Kolk: Lie Groups. Universitext,
Springer-Verlag 1999.

\bibitem{Scardua} B. Sc\'ardua: {\it Holomorphic Anosov Systems,
foliations and fiberings of complex manifolds}, Dynamicals
Systems, Vol. 18, No. 4, 2003, 365-389.


\bibitem{Steenrod} N.S. Steenrod: The Topology of fiber bundles,
Princeton Univ. Press, 1951.

\bibitem{Tischler} D. Tischler: {\it On fibering certain foliated
manifolds over $S^1$}, Topology, vol. 9 (1970), pp. 153-154.

\bibitem{V} V.S. Varadarajan: Lie Groups, Lie Algebras, and their
Representations. Prentice-Hall, Inc., Englewood Cliffs, N.J.,
1974.

\end{thebibliography}

\leftline{{\sc Alexandre  Behague and Bruno Sc\'ardua}}
\leftline{Instituto de Matem\'atica} \leftline{Universidade
Federal do Rio de Janeiro - Caixa Postal 68530}
\leftline{21945-970, Rio de Janeiro - RJ} \leftline{BRAZIL}
\leftline{abehague@ime.uerj.br and scardua@impa.br}

\end{document}